\documentclass[12pt]{amsart}
\newtheorem{thm}{Theorem}[section]
\newtheorem{tm}[thm]{Theorem}
\newtheorem{cor}[thm]{Corollary}
\newtheorem{lem}[thm]{Lemma}
\newtheorem{prop}[thm]{Proposition}
\theoremstyle{definition}

\theoremstyle{remark}

\numberwithin{equation}{section} % MATH

\def\r#1{\mathbb R^{#1}}

\def\X{\mathfrak X}
\def\id{\operatorname{id}}

\def\grad{\operatorname{grad}}

\begin{document}

\title{Transformations between surfaces in $R^4$ \\ with flat normal and/or tangent bundles}
\author{Angel Montesinos-Amilibia}
\address{Departamento de Geometria y Topologia, Universidad de Valencia, Vicente A. Estelles 1, 46100 Burjasot (Valencia), Spain}%
\email{montesin@uv.es}%
\date{\today}% %\dedicatory{}% %\commby{}% %
\thanks{Work partially supported by DGCYT grant n. BFM2003-02037.}
\begin{abstract} We exhibit several transformations of surfaces in
$R^4.$ First, one that takes a flat surface and gets a surface with
flat normal bundle; then, one that takes a surface with flat normal
bundle and gets a flat surface; finally, a one-parameter family of
transformations on a flat surface with flat normal bundle and gives
a flat surface with flat normal bundle. This family satisfies a
permutability property.
\end{abstract}

%----------------------------------------------------------------

\maketitle %
%----------------------------------------------------------------
\pagestyle{plain}

\section{Introduction}

Among surfaces in $\r4$, those of flat tangent bundle and those of
flat normal bundle have received considerable attention, especially
those that have both properties.

In this paper, we present several transformations between surfaces
of those types. In \cite{DoCarmo1993} we find  a transformation that
takes a hyper-spherical surface (hence, of flat normal bundle), and
gets a flat surface. We present here
 a transformation that takes any surface with flat normal bundle
without inflection points and converts it to its \it evolute, \rm
which results in a flat surface; the condition (no inflection
points) is meant to shun the possibility that the map go to
infinity, as happens in a point with zero curvature when defining
the evolute of a plane curve.

Then there is a kind of inverse, that is a transformation that takes
any immersed flat surface in $\r4$ and gives (in the region where
that transformation is an immersion) its \it envelope, \rm which is
a surface with flat normal bundle without inflection points.

Thus, it seems that the differential equations that define surfaces
in $\r4$ with flat normal bundle are essentially the same as those
that define flat surfaces.

By combining both types of transformations we get a transformation,
$f_t:M\to f_t(M),$ which depends on a real parameter $t,$ and yields
a flat surface with flat normal bundle without inflection points
from a surface with the same properties. These transformations
satisfy an analogous to the Bianchi permutability theorem for
B\"{a}cklund transformations (see \cite{Burs2004} for a description
in a broad context).

All the transformations $f:M\to f(M)$ so far described are
``orthogonal" in the sense that the tangent plane of $f(M)$ at
$f(p)$ is the orthogonal complement of the tangent plane of $M$ at
$p.$ The composition of two such transformations gives a ``parallel"
transformation, that is one such that the tangent plane of $f(M)$ at
$f(p)$ is equal to the tangent plane of $M$ at $p.$

\section{Basic concepts and notation}
In the following, $M$ will be a surface immersed in $\r n.$ However,
since all of our study will be local, one can without loss of rigor
assume that $M$ is an embedded surface. On $M$ we have the tangent
bundle $\pi:TM\to M,$ and the \it normal bundle \rm given by
$$NM= \cup_{p\in M}(T_pM)^\bot,\quad \pi_N: NM\to M,$$
where $(T_pM)^\bot$ denotes the subspace of $T_p\r n$ orthogonal to
$T_pM.$ Its fibre upon $p\in M$ will be denoted by
$N_pM=(T_pM)^\bot.$ Usually we will consider $T_pM$ and $N_pM$ as
vector subspaces of $\r n.$ We will use a dot to mean the standard
inner product. If $X\in T_p\r n,$ we will have $X= X^\top +X^\bot,$
with $X^\top\in T_pM,\; X^\bot\in N_pM.$

The Lie algebra of vector fields on a manifold $M$ will be denoted
$\X(M),$ and if $E$ is the total space of a vector bundle over $M$,
$\Gamma E$ will stand for the $C^\infty(M)$-module of its
differentiable sections. Usually, if $s$ is a section of a fiber
bundle, $s_p$ will be its value at $p.$

The ordinary directional derivative in $\r n$ will be written as
$D_X.$ But note that it may have a broader meaning of which we will
have a frequent use. In fact, if $S$ is a submanifold of $\r n,\;
p\in S,\; X_p\in T_pS$ and $f: S\to \r m$ is a differentiable map,
then $D_{X_p}f\in \r m$  will be defined as $df(X_p)\in T_{f(p)}\r
m\approx \r m.$ For vector fields on $\r m,\; D$ is a metric linear
connection with zero torsion and curvature.

Let $u\in \Gamma NM$ and $p\in M.$  The map $A_{u_p}: T_pM\to T_pM$
defined by $A_{u_p}(X)= (D_Xu)^\top,$ depends only on the value
$u_p$ of $u$ at $p$ and it is self-adjoint with respect to the inner
product. The second fundamental form of $M,\;\alpha,$ may be defined
at $p$ as the symmetric bilinear form $\alpha_p: T_pM\times T_pM\to
N_pM$ that satisfies $u_p\cdot \alpha_p(X_p,Y_p)= -A_{u_p}(X_p)\cdot
Y_p,$ for any $X_p,Y_p\in T_pM,\; u_p\in N_pM.$ We have equivalently
$$\alpha(X,Y)= (D_XY)^\bot,\quad X,Y\in \X(M).$$

The map $\X(M)\times\X(M)\to \X(M)$ given by
$\nabla^\top_XY=(D_XY)^\top$ is a torsionless metric linear
connection. If its curvature is zero everywhere, we say that $M$ is
flat; in this case, if $p$ is any point of $M,$ there is a chart
$\psi: (u,v)\in U\subset\r2\to M$ such that $p\in \psi(U)$ and that
the pull-back of the first fundamental form of $M$ by $\psi$ reads
$du^2+dv^2.$ The map $\X(M)\times\Gamma NM\to \Gamma NM$ given by
$\nabla^\bot_Xu=(D_Xu)^\bot$ is a metric linear connection. It its
curvature is zero everywhere, we say that the normal bundle of $M$
is flat.

Let $PT_pM$ the projective space of vector lines of $T_pM.$ The
second fundamental form defines a map $\eta_p:PT_pM\to N_pM,$  by
$$\eta_p([t])=\eta_p(t)=\frac{\alpha_p(t,t)}{t\cdot t},\quad t\in T_pM\backslash\{0\}.$$

The image of $\eta_p$ is an ellipse in $N_pM,$  called the \it
curvature ellipse \rm at $p.$ If that ellipse lies in an affine line
(i.e. it degenerates to a segment or to a point), we say that $p$ is
a \it semiumbilic \rm point. If, in addition, that line passes by
the origin of $N_pM,$ we say that $p$ is a \it point of inflection.
\rm

The following facts are well known or easily proved (see for
instance \cite{Nuno1996}, \cite{RFSB2002} and \cite{Mello2005}). The
point $p$ is semiumbilic iff the curvature of $\nabla^\bot$ vanishes
at $p.$ So, $M$ is totally semiumbilic iff its normal bundle is
flat. If $p$ is semiumbilic and $[t]\in PT_pM$ is such that
$\eta_p([t])$ is equal to any of the ends of the segment in which
the curvature ellipse degenerates, then we say that $[t]$ (or $t$)
is an \it asymptotic \rm direction at $p;$ if that segment does not
degenerate to a point at $p$ (\it umbilic \rm point) then there are
two asymptotic directions at $p$ and they are mutually orthogonal.
If $p$ is semiumbilic, then the curvature of $\nabla^\top$ vanishes
at $p$ iff the circle that has the curvature segment as diameter
passes by the origin of $N_pM;$ or, equivalently if p is not
umbilic, iff $\eta_p(t_1)$ and $\eta_p(t_2)$ are orthogonal, where
$t_1$ and $t_2$ are the asymptotic directions at $p.$ If, in
addition, $p$ is not of inflection, $\eta_p(t_1)$ and $\eta_p(t_2)$
are linearly independent.

\section{Some facts on semiumbilical surfaces}

The following characterization of semiumbilic points will be crucial
for our results. This definition differs from the first one that I
know, that of Wong \cite{W1945}. The reason is that for surfaces in
$\r n$ with $n\ge 5,$ Wong condition of being semiumbilic is
satisfied ``almost everywhere", that is, in addition to the points
that are semiumbilic for us, whenever the curvature ellipse does not
degenerate and the affine plane containing it does not pass by the
origin.

\begin{prop}
Let $M$ be an immersed surface in $\r n,\;n\ge 4,\; p\in M$ be a non
umbilic point and $g$ denote the first fundamental form of $M.$ Then
there is an affine plane $E_p$ that passes by the origin of $N_pM$
and contains the curvature ellipse at $p$ and there is a vector
$c_p\in E_p$ such that $c_p\cdot\alpha_p= g_p,$ iff $p$ is a non
inflection semiumbilic point. Moreover if such a vector $c_p$
exists, it is unique.
\end{prop}
\begin{proof}
If such a plane $E_p$ and vector $c_p$ exist, then for each unit
vector $t\in T_pM$ we will have
$c_p\cdot\alpha_p(t,t)=c_p\cdot\eta_p(t)= g_p(t,t)=1.$ Therefore,
the height of all points of the curvature ellipse at $p$ over the
hyperplane $H_p$ of $N_pM$ with normal $c_p$ is constant and equal
to $\frac1{|c_p|}.$ Therefore, the curvature ellipse lies in the
line of intersection of $E_p$ with the affine hyperplane parallel to
$H_p$ at that distance, so that $p$ is semiumbilic and obviously it
cannot be of inflection without being umbilic.

Conversely, if $p$ is a non inflection semiumbilic point, let $n_p$
be the point of the line containing the curvature ellipse at $p$ (a
segment) nearest to the origin of $N_pM.$ Since $p$ is not an
inflection point $n_p\ne 0.$ Let $c_p=\frac{n_p}{n_p\cdot n_p}.$
Since $n_p\cdot \eta_p(t)= n_p\cdot n_p,$ for all unit vectors $t\in
T_pM,$ se will have $c_p\cdot \eta_p(t)=1= g_p(t,t).$ Since
$c_p\cdot \eta_p$ is a quadratic form on $T_pM$, we conclude that
its corresponding bilinear symmetric form $c_p\cdot\alpha_p$ is
equal to $g_p.$ The uniqueness of $c_p$ is obvious.
\end{proof}

In $\r4$ we can dispense with the requirement of a plane $E_p$ to
contain $c_p$ because $N_pM$ does the job. Thus

\begin{cor}
Let $M$ be an immersed surface in $\r 4,\; p\in M$ be a non umbilic
point and $g$ denote the first fundamental form of $M.$ Then there
is a unique vector $c_p\in N_pM$ such that $c_p\cdot\alpha_p= g_p$
iff $p$ is a non inflection semiumbilic point.
\end{cor}

From now on, unless otherwise stated, $M$ will be a surface immersed
in $\r4.$ We describe now the curvature ellipse in more concrete
terms. If  $(t_1, t_2)$ is a local orthonormal frame of $TM,$ we put
$b_1=\eta(t_1)=\alpha(t_1,t_1),\; b_2=\eta(t_2)=\alpha(t_2, t_2), \;
b_3= \alpha(t_1,t_2).$ If $t\in \X(M)$ is a unit vector field, we
will have $t= t_1\cos\theta+ t_2\sin\theta.$ Then, we have
$\eta(t)=b_1\cos^2\theta+ b_2\sin^2\theta+ b_3\sin 2\theta.$ After
an easy calculation we get $$\eta(t)= H+ B\cos2\theta+
C\sin2\theta,$$ where
$$H= \frac12(b_1+b_2),\quad B=\frac12(b_1-b_2),\quad C= b_3$$
are smooth local sections of $NM$.

$H$ is called \it mean curvature vector (field) \rm and it does not
depend on the choice of the orthonormal frame $(t_1,t_2).$ The other
two sections $B$ and $C$ do depend on it. In a region where the
ellipse does not degenerate to a point or a circle, the frame $(t_1,
t_2)$ can be locally chosen so that the major half-axis of the
ellipse be $B$ and the minor, $C.$ That is  $|B|\ge|C|,$ and $B\cdot
C=0.$  We denote by $J: NM\to NM$ a local fibred map such that
$J(N_pM)=N_pM,\; J^2=-\id,\; Ju\cdot u=0,\; Ju\cdot Ju= u\cdot u,\;
\forall u\in NM.$ Then $J$ is defined up to a sign.

For avoiding annoying repetitions we will say that the immersed
surface $M$ is \it semiumbilical \rm if all of its points are non
inflection semiumbilic points.

\begin{prop} Let $M$ be an immersed semiumbilical surface in $\r4.$
Let the local orthonormal  frame $(t_1,t_2)$ of $TM$ satisfy
$b_3=\alpha(t_1,t_2)=0.$ Then the section
$$c=\frac{JB}{H\cdot JB},$$
where $B=\frac12(\alpha(t_1,t_1)-\alpha(t_2,t_2)),\;
H=\frac12(\alpha(t_1,t_1)+\alpha(t_2,t_2))$ is well defined and
satisfies $c\cdot\alpha=g.$\label{JB}
\end{prop}
\begin{proof}
Since there are no inflection points in $M,$ the curvature ellipse
at each point is a segment not collinear with the origin, that is
$b_1$ and $b_2$ are linearly independent. The frame $(t_1,t_2)$ is
simply a local frame of asymptotic directions. We have first $H\cdot
JB=\frac14(b_1+b_2)\cdot(Jb_1-Jb_2)=-\frac12 b_1\cdot Jb_2\ne 0$
because $b_1$ and $b_2$ are linearly independent. Therefore,
$$c\cdot\alpha(t_1,t_1)= \frac{Jb_1-Jb_2}{-b_1\cdot Jb_2}\cdot
b_1=1= t_1\cdot t_1.$$ In the same manner we get
$c\cdot\alpha(t_2,t_2)=t_2\cdot t_2=1.$ Since $\alpha(t_1,t_2)=0,$
we have finally $c\cdot\alpha = g.$\end{proof}

The following result will be used afterwards, in the context of
surfaces with both bundles, normal and tangent, flat.

\begin{lem} \label{2Lema} Let $M$ be a semiumbilical immersed
surface in $\r4.$ Let the local orthonormal frame $(t_1,t_2)$ of
 $TM$ satisfy $b_3=\alpha(t_1,t_2)=0$ and put $b_1=\alpha(t_1,t_1)$,
$b_2=\alpha(t_2,t_2).$ Then
$$b_1\cdot \nabla^\bot_{t_2}c = b_2\cdot \nabla^\bot_{t_1}c=0.$$
\end{lem}
\begin{proof}
We have $0=\alpha(t_1,t_2)= (D_{t_1}t_2)^\bot=(D_{t_2}t_1)^\bot.$
Therefore, $D_{t_1}t_2,\; D_{t_2}t_1\in\X(M).$ Now,
\begin{align*}
b_1\cdot \nabla^\bot_{t_2}c&=b_1\cdot D_{t_2}c=\alpha(t_1,t_1)\cdot
D_{t_2}c = (D_{t_1}t_1)^\bot\cdot D_{t_2}c\\
&= D_{t_2}((D_{t_1}t_1)^\bot\cdot c)-c\cdot
D_{t_2}(D_{t_1}t_1)^\bot\\
&=D_{t_2}(c\cdot\alpha(t_1,t_1))-c\cdot
D_{t_2}(D_{t_1}t_1-(D_{t_1}t_1)^\top)\\
&=D_{t_2}(1)- c\cdot D_{t_2}D_{t_1}t_1+
c\cdot\alpha(t_2,(D_{t_1}t_1)^\top)\\
&=-c\cdot(D_{t_1}D_{t_2}t_1+D_{[t_2,t_1]}t_1)+t_2\cdot D_{t_1}t_1.
\end{align*}
Now we observe that  $[t_2,t_1],\; D_{t_2}t_1\in \X(M),$ whence
$$b_1\cdot \nabla^\bot_{t_2}c=-t_1\cdot D_{t_2}t_1-[t_2,t_1]\cdot
t_1-t_1\cdot D_{t_1}t_2=0,$$ because $t_1\cdot D_{t_2}t_1=0$ and
$[t_2,t_1]= D_{t_2}t_1-D_{t_1}t_2.$ In the same manner we get
$b_2\cdot \nabla^\bot_{t_1}c=0.$
\end{proof}

\section{Transformations between surfaces with flat tangent bundle
and surfaces with flat normal bundle in $\r 4$}

A good part of our results are based in the following observations:

\begin{lem} Let $U,V$ be two surfaces immersed in $\r 4$ and let
$f:U\to V$ be a diffeomorphism such that for any  $p\in U$ we have
$T_pU= N_qV,$ where $q=f(p).$ Then, a section $Y$ of $TU$ satisfies
$\nabla^\top Y =0$ (is parallel) iff $\tilde Y= Y\circ f^{-1},$
which is a section of $NV,$ satisfies $\nabla^\bot\tilde Y=0$ (is
parallel). And a section $u$ of $NU$ is parallel iff $u\circ
f^{-1},$ which is a section of $TV,$ is parallel.\label{1lema}
\end{lem}
\begin{proof}
We have $\tilde Y_q= (Y\circ f^{-1})_{f(p)}= Y_p\in T_pU= N_qV;$
hence $\tilde Y$ is a section of $NV.$ Let $X\in T_pU= N_qV,\; u\in
N_pU= T_qV.$ We will have
$$X\cdot
\nabla^\bot_u\tilde Y=X\cdot D_u\tilde Y=X\cdot d\tilde Y(u)= X\cdot
dY( df^{-1}(u))= X\cdot
D_{df^{-1}(u)}Y=X\cdot\nabla^\top_{df^{-1}(u)}Y,$$ and now both
claims are evident.
\end{proof}

For the previous Lemma we did need that the ambient space were $\r
4$ because $$4= 2\dim(T_pU)=\dim(T_pU)+ \dim(N_qV)=\dim(T_pU)+
\dim(N_pU)$$ is the dimension of the ambient space. In the following
this condition is not necessary. The proof is similar.

\begin{lem} Let $U,V$ be two surfaces immersed in $\r n$ and let
$f:U\to V$ be a diffeomorphism such that for any  $p\in U$ we have
$T_pU= T_qV,$ where $q=f(p).$ Then, a section $Y$ of $TU$ is
parallel iff $\tilde Y= Y\circ f^{-1},$ which is a section of $TV,$
is parallel. And a section $u$ of $NU$ is parallel iff $u\circ
f^{-1},$ which is a section of $NV,$ is parallel.\label{3lema}
\end{lem}

Now we begin to study the conditions to have diffeomorphisms as
those used in the preceding Lemmas.

\begin{prop} Let $M$ be an immersed surface in $\r n.$ Then the
following statements are equivalent:
\begin{enumerate}
\item $M$ is flat;
\item Given any point $p\in M,$ there is an open neighborhood $U$
of $p$ in $M$ and a vector field $e\in\X(U)$ such that
$\nabla^\top_Xe= X$ for all $X\in \X(U);$
\item Given any point $p\in M,$ there is an open neighborhood $U$
of $p$ in $M$ and a vector field $e\in\X(U)$ such that for any $q\in
U$ we have $df(T_qU) \subset N_qU,$ where $f:U\to\r n$ is the map
defined by  $f(q)=q - e_q,$ that is $f=\id-e.$
\end{enumerate} \label{1prop}
\end{prop}
\begin{proof}
$(1)\Leftrightarrow(2).$ If $M$ is flat, there is a chart $\mu$ in
an open neighborhood $U$ of $p$, with coordinates $(u,v)$ whose
first fundamental form reads $du^2+dv^2$. This means that the vector
field $e= u\mu_u+ v\mu_v$ satisfies the required property because
$\mu_u, \mu_v$ are parallel. This field $e$ is determined upto the
addition of a parallel tangent vector field on $U$. Conversely, if
$e\in\X(U)$ satisfies the condition, then for any $X,Y\in\X(U)$ we
have
$$R^\top(X,Y)e=\nabla^\top_X\nabla^\top_Ye-
\nabla^\top_Y\nabla^\top_Xe-\nabla^\top_{[X,Y]}e=\nabla^\top_XY-
\nabla^\top_YX-[X,Y]=0,$$ because $\nabla^\top$ is torsionless.
Since $e$ can vanish only at isolated points and the dimension of
$M$ is $2,$ we conclude that $R^\top=0.$

$(2)\Leftrightarrow(3).$ Let $e:U\to TU$ be a local section of $TM$
and let $X,Y\in T_qM,\; q\in U.$ Then $Y\cdot df(X)= Y\cdot(X-
de(X))= Y\cdot(X-D_Xe)=Y\cdot(X-\nabla^\top_Xe)$ and now our claim
is evident.
\end{proof}

For $M$ flat, and assuming that the field $e$ is defined in all of
$M,$ let $f=\id-e:M\to \r n$ and assume that for some  $X\in T_pM$
we have $df(X)=0.$ This would be equivalent to say that for any
$u\in N_pM$ we had $0=u\cdot df(X)= u\cdot(X- D_Xe)= -u\cdot D_Xe=
-u\cdot\alpha_p(e_p,X)= 0,$ that is $\alpha_p(e_p,X) =0,$ that is
$\{0\}\ne\ker \alpha_p(e_p, ):T_pM\to N_pM.$ If $n=3,$ this always
happens, so that for that dimension $f$ is never an immersion.
Probably, if $n=4,\; df$ would generically not be one-to-one only
along some curves; if $n>4,\; f$ would be an immersion generically
outside isolated points, and so on. However, I shall not dwell on
this point.

Assume that in fact $f$ be an immersion, and let $\phi= (f|_U)^{-1}$
for some open $U\subset M$ such that $f|_U: U\to V=f(U)$ is a
diffeomorphism. Then $\phi=\id_V + e\circ\phi,$ as it can be proved
easily. Since $T_qV= df(T_pU)\subset N_pU,$ we will have
$T_pU\subset N_qV.$ Hence, $c=e\circ\phi$ is a section of $NV$ and
the map $\phi=\id+c$ satisfies $d\phi(T_qV)=T_pU\subset N_qV.$ This
motivates the following proposition.

\begin{prop}  Let $M$ be an immersed surface in $\r n$ with normal
bundle $NM$ and first and second fundamental forms $g$ and $\alpha,$
respectively. Let $c\in\Gamma NM$ and put $f=\id+c:M\to\r n.$ Then
$c\cdot\alpha = g$ iff for each $p\in M$ we have $df(T_pM) \subset
N_pM.$ \label{2prop}
\end{prop}
\begin{proof}
Let $p\in M,\; X,Y\in T_pM.$ Then
\begin{align*} Y\cdot
(df)_p(X)&=Y\cdot( X + (dc)_p(X))= Y\cdot(X + D_X c)\\
&=Y\cdot(X +A_c(X))= g(X,Y)-c\cdot\alpha(X,Y),
\end{align*}
and our claim is now evident.
\end{proof}

Assume now that $c\in\Gamma NM$ satisfies $c\cdot\alpha = g,$ and
let $f=\id+c:M\to\r n.$ If $p\in M,$ let us study the condition for
$df_p$ not being one-to-one. This happens iff there is some non
vanishing vector $X\in T_pM$ such that, for all $u\in N_pM$ the
following holds: $u\cdot df_p(X)=
u\cdot(X+D_Xc)=u\cdot\nabla^\bot_Xc=0.$ That is iff $\{0\}\ne \ker\;
(\nabla^\bot c)_p: X\in T_pM\mapsto \nabla^\bot_X c\in N_pM.$ As
before, we see that if $n=3,\; f$ cannot be an immersion, and that,
probably, for $n=4$ it fails generically to be so only on some
curves, etc.

The condition $c_p\cdot\alpha_p = g_p$ says that the height of the
curvature ellipse with respect to the vector hyperplane of $N_pM$
orthogonal to $c_p$ is constant and equal to $\frac1{|c_p|}$. If
$n=3$ this entails that $p$ is umbilic and not planar. If $n=4$ this
happens only if the ellipse degenerates to a point, not the origin,
or to an affine segment not collinear with the origin. If, finally,
$n > 4,$ this may occur almost always, because it is equivalent to
require only that the least affine subspace of $N_pM$ that contains
the curvature ellipse does not pass by the origin.

The next two Theorems, that are part of our main results explain why
from now on we consider only surfaces in $\r4.$ Roughly, they
establish a transformation of a surface with flat normal bundle to a
flat surface, and a transformation that takes a flat surface and
converts it to a surface with flat normal bundle.

\begin{tm}  Let $M$ be a semiumbilical surface immersed in $\r4.$
If $c$ is the section of $NM$ described in \ref{JB} that satisfies
$c\cdot\alpha=g,$ and $f=\id+c:M\to \r4$ is an immersion, then
$f(M)$ is an immersed flat surface that we call the \emph{evolute}
of $M.$  If, in addition, $M$ is flat, then $f(M)$ is semiumbilical.
\end{tm}
\begin{proof}
Let $U$ be an open subset of $M$ for which $f:U\to V=f(U)$ is a
diffeomorphism. We will have $\phi=(f|_U)^{-1}=\id_V-c\circ\phi,$
and  $d\phi(T_qV)\subset N_qV.$ Taking account of the dimensions, we
conclude that $df(T_pU)= N_pU= T_qV$ and $d\phi(T_qV)= T_pU= N_qV.$
By the same reason, $\tilde e=c\cdot\phi$ is a section of $TV$. Our
claims are now a consequence of \ref{1prop} and \ref{1lema}. In
fact, $M$ is flat iff there is a nonvanishing parallel vector field
on $M$, and all its points are semiumbilic iff its normal bundle is
flat, that is iff it admits a nonvanishing parallel section. In both
cases, due to the dimension 2 of those bundles. The question whether
$f(M)$ has inflection points when $M$ is semiumbilical and flat will
be settled below (Theorem \ref{2tm}(1)) under a more general
context.
\end{proof}

Had we assumed that $n>4,$ we could not conclude that $e$ was
tangent to $V$ nor that  $V$ was flat.

\begin{tm} Let $M$ be a surface immersed in $\r4.$ Let $e\in \X(M)$
be such that $\nabla^\top_Xe=X,\;\forall X\in\X(M)$ (hence, $M$ is
flat) and put $f=\id-e:M\to\r4.$ Then, if $f$ is an immersion, the
immersed surface $f(M)$ is semiumbilical and we will say that it is
an \emph{envelope} of $M.$ If, in addition, $M$ is semiumbilical,
then $f(M)$ is flat.
\end{tm}
\begin{proof}
If we put locally $\phi=f^{-1},$ and $c=e\circ\phi,$ then
$\phi=\id+c,$ and $d\phi(T_qf(M))= N_qf(M),$ and we conclude with
the aid of \ref{2prop} and \ref{1lema}.
\end{proof}

The evolute of a non-flat semiumbilical surface is unique. On the
contrary, there are many envelopes of a flat surface. In fact, it is
enough to add a parallel vector field $Z\in\X(M)$ to $e.$ However,
if a semiumbilical surface is also flat, then there may be more
evolutes. This is a consequence of this proposition that is slightly
more general.

\begin{prop}
Let $M$ be a flat semiumbilical surface in $\r n$, let $Z\in\X(M)$
be a parallel vector field and put $f=\id +c+Z:M\to \r n.$ Then
$df(T_pM)\subset N_pM$ for any $p\in M.$
\end{prop}
\begin{proof}
Let $X,Y\in T_pM.$ Using \ref{2prop}, we have:
$$Y\cdot df(X)= Y\cdot D_XZ= Y\cdot\nabla^\top_XZ=0.$$
\end{proof}

\section {Transformations of surfaces in $\r4$ with tangent and
normal bundles both flat}

As we have recalled, if a surface is flat and semiumbilical, there
is, in a neighborhood of any point, an orthonormal frame $(t_1,t_2)$
of asymptotic directions, such that, if
$b_1=\alpha(t_1,t_1)=(D_{t_1}t_1)^\bot,\;
b_2=\alpha(t_2,t_2)=(D_{t_2}t_2)^\bot,\;
b_3=\alpha(t_1,t_2)=(D_{t_1}t_2)^\bot,$ we have: \begin{enumerate}
\item $b_1$ and $b_2$ are linearly independent at each point;
\item $b_1\cdot b_2=0,\; b_3=0.$
\end{enumerate}

With this notation, we have:

\begin{prop} \label{4prop}
Let $M$ be a flat semiumbilical surface in $\r 4$, and let $c$ be
the section of $NM$ such that $c\cdot\alpha=g.$ Then,
\begin{enumerate}
 \item There is a unique vector field $j\in \X(M)$ such
that for any $X\in\X(M)$ we have $\alpha(j,X)=\nabla^\bot_Xc,$ and
it is given by
$$j=\frac12\grad(c\cdot c)=\frac{b_1\cdot D_{t_1}c}{b_1\cdot
b_1}t_1+\frac{b_2\cdot D_{t_2}c}{b_2\cdot b_2}t_2.$$
\item Let $f=\id + c+ Z:M\to \r4,$ where $Z\in\X(M)$ is parallel, and put $U$ to denote the open
subset of $M$ where $j+Z$ neither vanishes nor is parallel to an
asymptotic direction. Then, $f$ is an immersion on $U.$
\end{enumerate}

\end{prop}
\begin{proof}
(1) If such a vector field $j$ exists, then
\begin{align*}
j\cdot X&= c\cdot\alpha(j,X)=c\cdot\nabla^\bot_Xc=c\cdot
D_Xc=\frac12D_X(c\cdot c)=\frac12 d(c\cdot c)(X)\\
&=\frac12\grad(c\cdot c)\cdot X.
\end{align*}
Therefore, if it exists, $j$ is unique and is given by
$\frac12\grad(c\cdot c)$. Its existence is not evident. We can write
$j=j^1t_1+j^2t_2$ and must have
$\alpha(j,t_1)=j^1b_1=\nabla^\bot_{t_1}c.$ This reduces to the
following two conditions \begin{enumerate}
\item $j^1b_1\cdot b_1= b_1\cdot D_{t_1}c,\quad$ that is $\;\;j^1=
\frac{b_1\cdot D_{t_1}c}{b_1\cdot b_1}.$
\item $j^1b_1\cdot b_2=b_2\cdot D_{t_1}c.$
\end{enumerate}
Condition (1) determines $j^1.$ Since $b_1\cdot b_2=0,$ condition
(2) can be met iff $ b_2\cdot D_{t_1}c=0,$ but this is true by
\ref{2Lema}. The same happens to $j^2,$ and this proves our claims.
Note that we have found a vector field $j\in\X(M)$ such that
$$(D_X(c-j))^\bot=0,\;\forall X\in \X(M).$$

(2) Let $p\in M,\; 0\ne X\in T_pM,$ and assume that $df(X)=0.$ Since
$df(X)\in N_pM$, this happens iff for any $u\in N_pM$ we have
$u\cdot df(X)=0.$ But
$$u\cdot df(X)=u\cdot(X+D_Xc+ D_XZ)=u\cdot (D_X(c+Z))^\bot=u\cdot (D_X(j+Z))^\bot.$$
Hence, $df(X)=0$ iff $(D_X(j+Z))^\bot=0.$ Now, if $Z=
Z^1t_1+Z^2t_2,$ we have
\begin{align*}
(D_X(j+Z))^\bot&=((j^1+Z^1)(D_Xt_1)^\bot+(j^2+Z^2)(D_Xt_2)^\bot\\
&=(j^1+Z^1)X^1 b_1+(j^2+Z^2)X^2b_2.\end{align*} Since $b_1$ and
$b_2$ are linearly independent, this is zero iff $j^1+Z^1=0$ or
$j^2+Z^2 =0$ because $X\ne 0,$ and this proves our claim.
\end{proof}

\begin{cor} \label{2cor}
Let $M$ be a flat semiumbilical surface in $\r 4.$ Let $c$ be the
section of $NM$ such that $c\cdot\alpha=g,$ let
$j=\frac12\grad(c\cdot c)\in\X(M),$ and let $Z\in\Gamma NM$ be
parallel. If $f:M\to\r4$ is defined by $f=\id + t(c-j)+ Z,$ with
$t\in\r{},$ then $df(T_pM)\subset T_pM,\;\forall p\in M.$ Hence, if
$f$ is an immersion, the surface $f(M)$ is flat and semiumbilical.
\end{cor}

\begin{prop} Let $M$ be a flat semiumbilical surface in
$\r 4,\; e\in\X(M)$ be such that $\nabla^\top_X e=X,\;\forall
X\in\X(M),$ and let $U\subset M$ be the open region where $e$
neither vanishes nor is an asymptotic direction. Then $f=\id-e$ is
an immersion in $U$ and if $p_0\in U,$ then there is a section $k$
of $NM$ in a neighborhood $V$ of $p_0$ such that $(D_X(k-e))^\bot=0$
on $V.$ Obviously, $k$ is determined up to the addition of a
parallel section of $NM.$\end{prop}
\begin{proof}
First, we know that $df(T_pM)\subset N_pM.$ Let $(t_1,t_2)$ be a
local orthonormal basis of $T_pM$ of asymptotic directions. Let
$0\ne X\in T_pM.$ Then $df_p(X)=0$ iff for any $u\in N_pM$ we have
$0=u\cdot df_p(X)= u\cdot(X-D_Xe)= -u\cdot D_Xe=
-u\cdot\alpha(e,X),$ that is iff
$\alpha(e,X)=e^1X^1b_1+e^2X^2b_2=0,$ where $X^i, e^i$ are the
components of $X$ and $e$ in the basis $(t_1,t_2)$ and $b_i$ have
its usual meaning. Since $(b_1,b_2)$is linearly independent, we
conclude that $e^1=0$ or $e^2=0,$ because $X\ne 0.$ Therefore, $e$
vanishes at $p$ or determines an asymptotic direction.

The immersed surface $f(U)$ is flat and semiumbilical, as we know
already. Let $V$ be an open subset of $f(U)$ where there is $\tilde
e\in \X(V)$ such that $\nabla^\top_u\tilde e=u,\;\forall u\in
\X(V).$ We put $k=-\tilde e\circ f\in\Gamma N\phi(V),$ where
$\phi=f^{-1}.$ Then $$u= \nabla^\top_u\tilde
e=-(D_u(k\circ\phi))^\top=-(D_{d\phi(u)}k)^\bot.$$ Now
$u=df(d\phi(u))=(d\phi(u)-D_{d\phi(u)}e)^\bot=-(D_{d\phi(u)}e)^\bot,$
and our claims follow.
\end{proof}

We can combine both transformations:

\begin{thm}\label{1tm}
Let $M$ be a flat semiumbilical surface immersed in $\r4$ and let
$Z\in\Gamma NM$ be parallel. If the map $f:M\to\r4$ is defined by
$f=\id + t_1(e-k)+ t_2(c-j)+Z,$ with $t_1, t_2\in\r{},$ then
$df(T_pM)\subset T_pM,\;\forall p\in M,$ and if $f: M\to f(M)$ is a
diffeomorphism, then $f(M)$ is flat and semiumbilical.
\end{thm}
\begin{proof}
Let $u\in N_pM,\; X\in T_pM.$ We will have \begin{align*}u\cdot
&df(X)= u\cdot(X+ t_1 D_X(e-k)+t_2D_X(c-j)+
D_XZ)\\
&=u\cdot((t_1 D_X(e-k))^\bot +t_2(D_X(c-j))^\bot+\nabla^\bot_XZ)=0.
\end{align*}
Now $f(M)$ is flat and with flat normal bundle as a consequence of
\ref{3lema}. We need still to prove that there are no inflection
points in $f(M).$ Let $p\in M,\; q=f(p),\; X,Y\in\X(M),\; u\in
\Gamma NM.$ Then
$(u\circ\phi)\cdot\tilde\alpha(X\circ\phi,Y\circ\phi)=
(u\circ\phi)\cdot D_{X\circ\phi}(Y\circ\phi).$ The evaluation of
this at $q$ gives $$u_p\cdot\tilde\alpha(X_p,Y_p)=u_p\cdot
D_{X_p}(Y\circ\phi)=u_p\cdot
D_{d\phi(X_p)}Y=u_p\cdot\alpha(d\phi(X_p),Y_p).$$ Therefore
$\tilde\alpha(X_p,Y_p)=\alpha(d\phi(X_p),Y_p).$ The curvature
ellipse at $q$ lies in an affine line of $N_qf(M).$ If $q$ is an
inflection point of $f(M)$ (this includes the case that it be
umbilic) then that line passes by the origin (the ellipse itself
passes by it because $f(M)$ is flat). Hence, the subspace
$\tilde\alpha(T_qf(M)\otimes T_qf(M))\subset N_qf(M)$ is
one-dimensional. The same must be true for $\alpha,$ so that $p$
would be an inflection point, against our hypotheses. Hence, $f(M)$
is semiumbilical.
\end{proof}

This transformation is said to be of \it parallel \rm type.

Now, we exhibit a transformation that sends each tangent space to
its orthogonal: this is another of our main results. Let us denote
with the same letter, crowned by a tilde, functions on $f(M)$ that
correspond to functions on $M.$

\begin{tm} \label{2tm}  Let $M$ be a flat semiumbilical surface
in $\r 4.$ Let $c$ be the section of $NM$ such that
$c\cdot\alpha=g,\; e\in\X(M)$ such that $\nabla^\top_Xe=X,\;\forall
X\in\X(M),$  and $t\in\r{}.$ Put $f=\id+tc-(1-t)e:M\to\r4.$ Then
$df(T_pM)\subset N_pM$ for any $p\in M.$ If $f$ is an immersion,
then the surface $f(M)$ is flat and semiumbilical, and if $\phi$
denotes a local inverse of $f$, we have
\begin{enumerate}
\item  If $X_p\in N_{f(p)}f(M)=T_pM,\; u_p,v_p\in
T_{f(p)}f(M)=N_pM,$ then
\begin{align*}
X_p\cdot\tilde\alpha(u_p,v_p)&=-v_p\cdot\alpha(X_p, d\phi(u_p)),\\
 X_p\cdot\tilde\alpha(df(Y_p),v_p)&=-v_p\cdot\alpha(X_p,Y_p).
 \end{align*}
 \item $\tilde c=((1-t)e -\frac t2\grad(c\cdot c))\circ\phi=((1-t)e -tj)\circ\phi.$
\item If $t\ne 1,$ there is a section $k\in\Gamma NM$ such that $(D_X(e-k))^\bot=0.$
This section is determined under the addition of a parallel section
of $NM.$ We have $\tilde e= (tc-(1-t)k)\circ\phi.$
\end{enumerate}
\end{tm}
\begin{proof}
Let $X,Y\in T_pM.$ Then: \begin{align*}X\cdot df(Y)&=X\cdot (Y+
t dc(Y)-(1-t) de(Y))\\
&=X\cdot(Y+tD_Yc-(1-t)D_Ye)\\
&=X\cdot Y- tc\cdot\alpha(X,Y)-(1-t)X\cdot\nabla^\top_Y
e\\
&=X\cdot Y(1-t-(1-t))=0,
\end{align*}
so that $df(T_pM)\subset N_pM,$ as claimed. If $f$ is an immersion,
then $f(M)$ is flat and with flat normal bundle. Now

(1) Let $u,v\in\Gamma(NM)$ and $X\in\X(M).$ Then $u\circ\phi,
v\circ\phi\in \X(f(M))$ and $X\circ\phi\in\Gamma(Nf(M)).$ Thus
\begin{align*}
(X\circ\phi)&\cdot\tilde\alpha(u\circ\phi, v\circ\phi)=
-(v\circ\phi)\cdot D_{u\circ\phi}(X\circ\phi)=-(v\circ\phi)\cdot
d(X\circ\phi)(u\circ\phi)\\
&=-(v\circ\phi)\cdot D_{d\phi\circ
u\circ\phi}X=-(v\circ\phi)\cdot\alpha(X\circ\phi,d\phi\circ
u\circ\phi)
\end{align*}
If $p\in M$ y $q=f(p),$ the value of the above expression at $q$
reads
$$X_p\cdot\tilde\alpha(u_p,v_p)=-v_p\cdot\alpha(X_p,
d\phi(u_p)),$$ or also, putting $u_p= df(Y_p)$:
$$X_p\cdot\tilde\alpha(df(Y_p),v_p)=-v_p\cdot\alpha(X_p,Y_p).$$

If $q\in f(M)$ then the curvature ellipse at $q$ lies in an affine
line of $N_qf(M).$ If $q$ is an inflection point of $f(M)$ (this
includes the case that it be umbilic) then that line passes by the
origin (the ellipse itself passes by it because $f(M)$ is flat).
Hence, $\tilde\alpha(df(Y_p),v_p)$ lies in the vector line
containing the ellipse for any $Y_p,v_p.$ If $X_p$ is non zero and
orthogonal to that line, then $\alpha(X_p,Y_p)=0$ for all $Y_p.$ If
we choose an orthonormal basis $(t_1,t_2)$ of $T_pM$ that determines
the asymptotic directions at that point, and $X_p= X^1 t_1+X^2t_2,\;
Y_p= Y^1 t_1+Y^2 t_2,$ then $\alpha(X_p,Y_p)=X^1Y^1b_1+X^2Y^2b_2=0.$
Since $(b_1,b_2)$ is linearly independent, we have
$X^1Y^1=X^2Y^2=0,$ and being $Y_p$ arbitrary we conclude $X_p=0,$
against our hypotheses. Hence, $f(M)$ is semiumbilical.

 (2) By the above result, there is $\tilde c\in \Gamma Nf(M)$ such that
$\tilde c\cdot \tilde \alpha=\tilde g.$ Assume in the above formula
that $X_p= \tilde c_q.$ Then $df(Y_p)\cdot v_p=
-v_p\cdot\alpha(\tilde c_q,Y_p),$ that is
\begin{align*} \alpha(\tilde c_q,Y_p)&= -df(Y_p)= -(Y_p
+tD_{Y_p}c-(1-t)D_{Y_p}e)^\bot\\
&=-t \nabla^\bot_{Y_p}c +(1-t)\alpha(Y_p,e_p), \end{align*} or also
$$\alpha(\tilde c_q-(1-t)e_p,Y_p)=-t \nabla^\bot_{Y_p}c=-\alpha(tj,Y_p)$$
by Proposition \ref{4prop}. So, $\tilde c_q=(1-t)e_p -tj_p,$ that is
$$\tilde c=((1-t)e -\frac t2\grad(c\cdot c))\circ\phi=((1-t)e -tj)\circ\phi.$$

(3) Since $f(M)$ is flat, there is $\tilde e\in \X(f(M))$ such that
$\nabla^\top_u \tilde e=u$ for any $u\in \X(f(M)).$ Let us define
$k\in \Gamma NM$ by means of $(1-t)k\circ \phi=tc\circ\phi - \tilde
e.$ Using the same technique we have on one hand:
\begin{align*} \alpha(\tilde c_q,Y_p)&= -df(Y_p)=-(D_{df(Y_p)}\tilde e)^\bot
=-(D_{Y_p}(\tilde e\circ f))^\bot
\\&=-t(D_{Y_p}c)^\bot + (1-t)(D_{Y_p}k)^\bot
\end{align*}
and on the other
$$\alpha(\tilde c_q,Y_p)=-t \nabla^\bot_{Y_p}c
+(1-t)\alpha(Y_p,e_p).$$ That is $k$ satisfies
$$(D_{Y_p}k)^\bot-\alpha(Y_p,e_p)=(D_{Y_p}(k-e))^\bot=0,$$
as desired.
\end{proof}

Now we will see that the composition of two transformations of
orthogonal type is one of parallel type, as we could expect. Let
$f_1=\id + t_1c-(1-t_1)e:M\to\r4$ and assume that $f_1$ is a
diffeomorphism. Since then $f(M)$ is flat and semiumbilical we can
transform it by $f_2=\id_{f_1(M)}+t_2\tilde c-(1-t_2)\tilde e,$ for
obtaining another flat and semiumbilical surface. We have, for $p\in
M, q=f(p)$:
\begin{align*}
(f_2\circ & f_1)(p)= f_2(p+ t_1c_p-(1-t_1)e_p)= p+
t_1c_p-(1-t_1)e_p+t_2\tilde c_q-(1-t_2)\tilde e_q\\
&=p+ t_1c_p-(1-t_1)e_p+t_2((1-t_1)e_p-t_1
j_p)-(1-t_2)(t_1c_p-(1-t_1)k_p)\\
&=p + t_1t_2(c_p-j_p)-(1-t_1)(1-t_2)(e_p-k_p)= (f_1\circ f_2)(p),
\end{align*}
where the final equality must be taken modulo the addition of a
parallel section of $NM.$   We have arrived thus to a result similar
to the Bianchi permutability theorem for B\"{a}cklund
transformations (see \cite{Burs2004})

\end{document}